\documentclass[12pt,a4paper]{amsart}

\usepackage{amscd}
\usepackage{graphicx}
\usepackage{amssymb,amsfonts,latexsym}
\usepackage{geometry}

\input xy
\xyoption {all}

\newtheorem{thm}{Theorem}[section]
\newtheorem{cor}[thm]{Corollary}
\newtheorem{lem}[thm]{Lemma}
\newtheorem{prop}[thm]{Proposition}
\theoremstyle{definition}

\newtheorem{defn}[thm]{Definition}
\theoremstyle{remark}
\newtheorem{rem}[thm]{Remark}

\newtheorem{exam}[thm]{Example}
\numberwithin{equation}{section}

\newcommand{\Q}{\mathbb{Q}}
\newcommand{\mQ}{\mathbb{Q}}
\newcommand{\ra}{\rightarrow}

\newcommand{\cO}{\mathcal{O}}
\newcommand{\cS}{\mathcal{S}}
\newcommand{\cA}{\mathcal{A}}
\newcommand{\pht}{\widetilde{\phi}}
\newcommand{\phat}{\hat{\phi}}
\newcommand{\pst}{\widetilde{\psi}}
\newcommand{\psat}{\hat{\psi}}
\newcommand{\fp}{\frak{p}}

\newcommand{\fa}{\frak{a}}

\long\def\forget#1\forgotten{}

\DeclareMathOperator\dl{dl}
\DeclareMathOperator\wl{wl}
\DeclareMathOperator\dg{d}
\DeclareMathOperator\kernel{ker}

\begin{document}

\title[minimal ramification problem]{On the minimal ramification problem for semiabelian groups}
\thanks{The research of the first author was supported in part by a grant from NSERC}

%
%
%
%
%
%

\def\Tech{Department of Mathematics, Technion-Israel Institute of Technology, Haifa 32000, Israel}
\def\Conc{Department of Mathematics and Statistics, Concordia
University, Montreal, Quebec, Canada}
\author{ Hershy Kisilevsky}
\address{\Conc}
\email{kisilev@mathstat.concordia.ca}

\author{ Danny Neftin}
\address{\Tech}
\email{neftind@tx.technion.ac.il}%

\author{ Jack Sonn }
\address{\Tech}
\email{sonn@math.technion.ac.il}

\begin{abstract}
It is known
(\cite{KS}, \cite{N})
that for any prime $p$ and any finite semiabelian $p$-group $G$, there exists a (tame) realization
of $G$ as a Galois group over the rationals $\Bbb Q$ with exactly $d=\dg(G)$ ramified primes, where
$\dg(G)$ is the minimal number of generators of $G$, which solves the minimal ramification problem
for finite semiabelian $p$-groups.  We generalize this result to obtain a theorem on finite
semiabelian groups and derive the solution to the minimal ramification problem for a certain family
of semiabelian groups that includes all finite nilpotent semiabelian groups $G$.  Finally, we give some indication of
the depth of the minimal ramification problem for semiabelian groups not covered by our theorem.
\end{abstract}

\maketitle

\section{Introduction}

Let $G$ be a finite group.  Let $d=\dg(G)$ be the smallest number for which there exists a subset $S$
of $G$ with $d$ elements such that the normal subgroup of $G$ generated by $S$ is all of $G$.  One
observes that if $G$ is realizable as a Galois group $G(K/\mQ)$ with $K/\mQ$ tamely ramified (e.g.
if none of the ramified primes divide the order of $G$), then at least $\dg(G)$ rational primes
ramify in $K$ (see e.g. \cite{KS}).  The \it minimal ramification problem for $G$ \rm is to realize
$G$ as the Galois group of a tamely ramified extension $K/\mQ$ in which exactly $\dg(G)$ rational
primes ramify. This variant of the inverse Galois problem is open even for $p$-groups, and no
counterexample has been found.  It is known that the problem has an affirmative solution for
all semiabelian $p$-groups, for all rational primes $p$ (\cite{KS},\cite{N}).  A finite
group $G$ is \it semiabelian \rm if and only if $G \in \cS\cA$, where $\cS\cA$ is the smallest
family of finite groups satisfying: (i) every finite abelian group belongs to $\cS\cA$. (ii) if $G
\in \cS\cA$ and $A$ is finite abelian, then any semidirect product $A \rtimes G$ belongs to
$\cS\cA$. (iii) if $G \in \cS\cA$, then every homomorphic image of $G$ belongs to $\cS\cA$. In this
paper we generalize this result to arbitrary finite semiabelian groups by means of a ``wreath
product length" $\wl(G)$ of a finite semiabelian group $G$.  When a finite semiabelian group $G$ is
nilpotent, $\wl(G)=\dg(G)$, which for nilpotent groups $G$ equals the (more familiar) minimal number
of generators of $G$.  Thus the general result does not solve the minimal ramification problem for
all finite semiabelian groups, but does specialize to an affirmative solution to the minimal
ramification problem for nilpotent semiabelian groups. Note that for a nilpotent group $G$, $\dg(G)$ is $max_{p||G|}\dg(G_p)$ and not $\sum_{p||G|}\dg(G_p)$, where $G_p$ is the $p$-Sylow subgroup of $G$. Thus, a solution to the minimal ramification problem for nilpotent groups does not follow trivially from the solution for $p$-groups.

\section{Properties of wreath products}

\subsection{Functoriality}
The family of semiabelian groups can also be defined using wreath
products. Let us recall the definition of a wreath product. Here and
throughout the text the actions of groups on sets are all right
actions. \begin{defn} Let $G$ and $H$ be two groups that act on the
sets $X$ and $Y$, respectively. The \it (permutational) wreath
product \rm $H\wr_X G$ is the set $H^X\times G=\{(f,g)|f:X\ra H,g\in
G\}$ which is a group with respect to the multiplication: $$
(f_1,g_1)(f_2,g_2)=(f_1f_2^{g_1^{-1}},g_1g_2),$$ where
$f_2^{g_1^{-1}}$ is defined by $f_2^{g_1^{-1}}(x)=f_2(xg_1)$ for any
$g_1,g_2\in G,x\in X,f_1,f_2:X\ra H$. The group $H\wr_X G$ acts on
the set $Y\times X$ by $(y,x)\cdot (f,g)=(yf(x),xg)$, for any $y\in
Y,x\in X,f:X\ra H,g\in G$. \end{defn} 
\begin{defn}The \textit{standard (or regular) wreath product} $H\wr G$ is defined
as the permutational wreath product with $X=G$,$Y=H$ and the right regular actions. \end{defn}

The functoriality of the arguments of a wreath product will play an important
role in the sequel. The following five lemmas are devoted to these functoriality properties.

\begin{defn} Let $G$ be a group that acts on $X$ and $Y$. A map
$\phi:X\ra Y$ is called a $G$-map if $\phi(xg)=\phi(x)g$ for every
$g\in G$ and $x\in X$. \end{defn} Note that for such $\phi$, we also
have
$\phi^{-1}(y)g=\{xg|\phi(x)=y\}=\{x'|\phi(x'g^{-1})=y\}=\{x'|\phi(x')=yg\}=\phi^{-1}(yg).$
\begin{lem}\label{G-map induces an epimophism} Let $G$ be a group
that acts on the finite sets $X,Y$ and let $A$ be an abelian group.
Then every $G$-map $\phi:X\ra Y$ induces a homomorphism $\pht:
A\wr_X G\ra A\wr_Y G$ by defining: $(\pht(f,g)) = (\hat{\phi}(f),g)$
for every $f:X\ra A$ and $g\in G$, where $\hat{\phi}(f):Y\ra A$ is
defined by: $$\hat{\phi}(f)(y)=\prod_{x\in \phi^{-1}(y)}f(x),$$ for
every $y\in Y$. Furthermore, if $\phi$ is surjective then $\pht$ is
an epimorphism. \end{lem} \begin{proof} Let us show the above $\pht$
is indeed a homomorphism. For this we claim:
$\pht((f_1,g_1)(f_2,g_2))=\pht(f_1,g_1)\pht(f_2,g_2)$ for every
$g_1,g_2\in G$ and $f_1,f_2:X\ra A$. By definition:
\begin{equation*}
\pht(f_1,g_1)\pht(f_2,g_2)=(\phat(f_1),g_1)(\phat(f_2),g_2)=(\phat(f_1)\phat(f_2)^{g_1^{-1}},g_1g_2),\end{equation*}
while:
$\pht((f_1,g_1)(f_2,g_2))=\pht(f_1f_2^{g_1^{-1}},g_1g_2)=(\phat(f_1f_2^{g_1^{-1}}),g_1g_2).$
We shall show that $\phat(f_1f_2)=\phat(f_1)\phat(f_2)$ and
$\phat(f^g)=\phat(f)^g$ for every $f_1,f_2,f:X \ra A$ and $g\in G$.
Clearly this will imply the claim. The first assertion follows
since: \begin{equation*}\phat(f_1f_2)(y)=\prod_{x\in \phi^{-1}(y)}
f_1(x)f_2(x)=\prod_{x\in \phi^{-1}(y)} f_1(x)\prod_{x\in
\phi^{-1}(y)} f_2(x)=\end{equation*}$$ =
\phat(f_1)(y)\phat(f_2)(y).$$ As to the second assertion we have:
\begin{equation*} \phat(f^g)(y)= \prod_{x\in \phi^{-1}(y)}
f^g(x)=\prod_{x\in \phi^{-1}(y)} f(xg^{-1})
=\end{equation*}\begin{equation*} = \prod_{x'g\in \phi^{-1}(y)}
f(x')  = \prod_{x'\in \phi^{-1}(y)g^{-1}} f(x').\end{equation*}
Since $\phi$ is a $G$-map we have
$\phi^{-1}(y)g^{-1}=\phi^{-1}(yg^{-1})$ and thus \begin{equation*}
\phat(f^g)(y)= \prod_{x\in \phi^{-1}(y)g^{-1}} f(x)= \prod_{x\in
\phi^{-1}(yg^{-1})} f(x) =  \phat(f)^g(y). \end{equation*} This
proves the second assertion and hence the claim. It is left to show
that if $\phi$ is surjective then $\pht$ is surjective. Let $f':Y\ra
A$ and $g'\in G$. Let us define an $f:X\ra A$ that will map to $f'$.
For every $y\in Y$ choose an element $x_y\in X$ for which
$\phi(x_y)=y$ and define $f(x_y):=f'(y)$. Define $f(x)=1$ for any
$x\not\in\{x_y|y\in Y\}$. Then clearly
$$\phat(f)(y)=\prod_{x\in\phi^{-1}(y)}f(x)=f(x_y)=f'(y).$$ Thus,
$\pht(f,g')=(\phat(f),g')=(f',g')$ and $\pht$ is onto. \end{proof}
\begin{lem}\label{construction of G-maps} Let $B$ and $C$ be two
groups. Then there is a surjective $B\wr C$-map
$\phi:B\wr C\ra B\times C$ defined by: $\phi(f,c)=(f(1),c)$ for
every $f:C\ra B, c\in C$. \end{lem} \begin{proof} Let
$(f,c),(f',c')$ be two elements of $B\wr C$. We check that
$\phi((f,c)(f',c'))=\phi(f,c)(f',c')$. Indeed, \begin{equation*}
\phi((f,c)(f',c'))=\phi(ff'^{c^{-1}},cc')=(f(1)f'^{c^{-1}}(1),cc')=(f(1)f'(c),cc')=\end{equation*}
$$ =(f(1),c)(f',c)=\phi(f,c)(f',c'). $$ Note that the map
$\phi$ is surjective: For every $b\in B$ and $c\in C$, one can
choose a function $f_b:C\ra B$ for which $f_b(1)=b$. One has:
$\phi(f_b,c)=(b,c)$. \end{proof}

The following Lemma appears in \cite[Part I, Chapter I, Theorem 4.13]{B} and describes
the functoriality of the first argument in the wreath product.
\begin{lem}\label{functoriality of
first argument} Let $G,A,B$ be groups and $h:A\ra B$ a homomorphism (resp. epimorphism).
Then there is a naturally induced homomorphism (resp. epimorphism) $h_*:A\wr G\ra B\wr
G$ given by $h_*(f,g)=(h\circ f,g)$ for every $g\in G$ and $f:G\ra
A$. \end{lem} 
%

The functoriality of the second argument is given in \cite[Lemma 2.15]{N} whenever the first argument is abelian:
\begin{lem}\label{functoriality with
second argument} Let $A$ be an abelian group and let $\psi:G\ra H$
be a homomorphism (resp. epimorphism) of finite groups. Then there is a homomorphism (resp. epimorphism)
$\pst:A\wr G\ra A\wr H$ that is defined by:
$\pst(f,g)=(\psat(f),\psi(g))$ with $\psat(f)(h)=\prod_{k\in
\psi^{-1}(h)}f(k)$ for every $h\in H$. \end{lem} 

These functoriality properties can now be joined to give a connection between different bracketing of iterated wreath products:
\begin{lem}
\label{induction step}  Let $A,B,C$ be finite groups and $A$
abelian. Then there are epimorphisms: $$A\wr (B \wr C)\ra (A\wr
B)\wr C\ra (A\times B)\wr C.$$ \end{lem}
\begin{proof} Let us first construct an
epimorphism $h_*:(A\wr B)\wr C\ra (A\times B)\wr C$. Define $h:A\wr
B\ra A\times B$ by: $$h(f,b)=(\prod_{x\in B}f(x),b),$$ for any
$f:B\ra A,b\in B$. Since $A$ is abelian $h$ is a homomorphism. For
every $a\in A$, let $f_a:B\ra A$ be the map $f_a(b')=0$ for any
$1\not=b'\in B$ and $f_a(e)=a$.   Then clearly $h(f_a,b)=(a,b)$
for any $a\in A,b\in B$ and hence $h$ is onto. By Lemma
\ref{functoriality of first argument}, $h$ induces an epimorphism
$h_*:(A\wr B)\wr C\ra (A\times B)\wr C$. To construct the
epimorphism $A\wr (B \wr C)\ra (A\wr B)\wr C$, we shall use the
associativity of the permutational wreath product (see \cite[Theorem 3.2]{B}).
Using this associativity one has: $$(A\wr B)\wr C= (A\wr_B
B)\wr_C C\cong A\wr_{B\times C} (B\wr_C C).$$ It is now left to construct an
epimorphism: $$ A\wr (B\wr C)=A\wr_{B\wr C} (B\wr C)\ra
A\wr_{B\times C} (B\wr C).$$ By Lemma \ref{construction of G-maps},
there is a $B\wr C$-map $\phi:B\wr C\ra B\times C$ and hence by
Lemma \ref{G-map induces an epimophism} there is an epimorphism
$A\wr_{B\wr C} (B\wr C)\ra A\wr_{B\times C} (B\wr C)$. \end{proof}

Let us iterate Lemma \ref{induction step}. Let $G_1,...,G_n$ be groups. The \it ascending iterated standard wreath product \rm of $G_1,...,G_n$ is
defined as $$(\cdots ((G_1 \wr G_2) \wr G_3)\wr \cdots) \wr G_n,$$ and the  \it descending iterated
standard wreath product \rm of $G_1,...,G_n$ is defined as $$G_1 \wr (G_2 \wr (G_3 \wr \cdots \wr
G_n))\cdots ).$$  These two iterated wreath products are not isomorphic in general, as the standard wreath
product is not associative (as opposed to the ``permutation" wreath product). We shall abbreviate and write $G_1\wr (G_2\wr ... \wr G_n)$ to refer to the descending wreath product and $(G_1 \wr ... \wr G_{r-1})\wr G_r$ to refer to the ascending wreath product. 

By iterating the epimorphism in Lemma
\ref{induction step} one obtains:
\begin{cor}\label{decending to ascending} Let $A_1,..,A_r$ be abelian groups. Then $(A_1 \wr ... \wr A_{r-1})\wr A_r$ is an  epimorphic image of  $A_1 \wr ( A_2 \wr ... \wr A_r).$
\end{cor}
\begin{proof} By induction on $r$. The cases $r=1,2$ are trivial; assume $r\geq 3$. By the induction hypothesis there is an epimorphism $$\pi_1':A_1 \wr (A_2\wr .... \wr A_{r-1})\ra (A_1\wr ... \wr A_{r-2})\wr A_{r-1}.$$ By Lemma \ref{functoriality of first argument}, $\pi_1'$ induces an epimorphism $\pi_1:(A_1\wr (A_2\wr ... \wr A_{r-1}))\wr A_r \ra (A_1\wr ... \wr A_{r-1})\wr A_r$.
Applying Lemma \ref{induction step} with $A=A_1,B=A_2\wr (A_3\wr ... \wr A_{r-1}), C=A_r$, one obtains an epimorphism:
$$\pi_2: A_1\wr (A_2\wr ... \wr A_r)\ra (A_1\wr (A_2\wr ... \wr A_{r-1}))\wr A_r.$$ Taking the composition $\pi=\pi_1\pi_2$ one obtains an epimorphism $$\pi:A_1\wr (A_2\wr ... \wr A_r)\ra (A_1\wr ... \wr A_{r-1})\wr A_r.$$

\end{proof}
\subsection{Dimension under epimorphisms}
Let us understand how the ``dimension" $\dg$ behaves under the homomorphisms in Lemma \ref{induction step} and Corollary \ref{decending to ascending}. By \cite{KL}, for any finite group $G$ that is not perfect, i.e. $[G,G]\not=G$, where $[G,G]$ denotes the commutator subgroup of $G$, one has $\dg(G)=\dg(G/[G,G])$. According to our definitions, for a perfect group $G$, $\dg(G/[G,G])=\dg(\{1\})=0$, but if $G$ is nontrivial,  $\dg(G)\geq 1$. As nontrivial semiabelian groups are not perfect, this difference will not effect any of the arguments in the sequel.
\begin{defn} Let $G$ be  a finite group and $p$ a prime. Define $\dg_p(G)$ to be the rank of the $p$-Sylow subgroup of $G/[G,G]$, i.e. $\dg_p(G):=\dg((G/[G,G])(p))$.
\end{defn}
Note that if $G$ is not perfect one has $\dg(G)=max_{p}(\dg_p(G))$. 

Let $p$ be a prime. An epimorphism $f:G\ra H$ is called $\dg$-preserving (resp. $\dg_p$-preserving) if $\dg(G)=\dg(H)$ (resp. $\dg_p(G)=\dg_p(H)$).

\begin{lem} Let $G$ and $H$ be two finite groups. Then:
$$H\wr G/[H\wr G,H\wr G]\cong H/[H,H]\times G/[G,G].$$
\end{lem}

\begin{proof} Applying Lemmas \ref{functoriality of first argument} and \ref{functoriality with second argument} one obtains an epimorphism $$H\wr G\ra H/[H,H]\wr G/[G,G].$$ By Lemma \ref{induction step} (applied with $C=1$) there is an epimorphism $$H/[H,H]\wr G/[G,G]\ra H/[H,H]\times G/[G,G].$$ Composing these epimorphisms one obtains an epimorphism $$\pi:H\wr G\ra H/[H,H]\times G/[G,G],$$
that sends an element $(f:G\ra H,g)\in H\wr G$ to $$(\prod_{x\in G}f(x)[H,H],g[G,G])\in H/[H,H]\times G/[G,G].$$
The image of $\pi$ is abelian and hence $\kernel(\pi)$ contains $K:=[H\wr G,H\wr G]$.

Let us show $K\supseteq \kernel(\pi)$. Let $(f,g)\in \kernel(\pi)$. Then $g\in [G,G]$ and $\prod_{x\in G}f(x)\in [H,H]$. As $g\in [G,G]$, it suffices to show that the element $f=(f,1)\in H\wr G$ is in $K$. Let $g_1,...,g_n$ be the elements of $G$ and for every $i=1,...,n,$ let $f_i$ be the function for which $f_i(g_i)=f(g_i)$ and $f(g_j)=1$ for every $j\not= i$. One can write $f$ as $\prod_{i=1}^n f_i$. Now for every $i=1,..,n$, the function $f_{1,i}=f_i^{g_i^{-1}}$ satisfies $f_{1,i}(1)=f(g_i)$ and $f_{1,i}(g_j)=1$ for every $j\not=1$. Thus $f_i$ is a product of an element in $[H^{|G|},G]$ and $f_{i,1}$. So,  $f$ is a product of elements in $[H^{|G|},G]$ and $f'=\prod_{i=1}^n f_{1,i}$. But $f'(1)=\prod_{x\in G}f(x)\in [H,H]$ and $f'(g_i)=1$ for every $i\not= 1$ and hence $f'\in [H^{|G|},H^{|G|}]$. Thus, $f\in K$ as required and $K=\kernel{\pi}$.

\end{proof} The following is an immediate conclusion:
\begin{cor} Let $G$ and $H$ be two finite groups. Then $$\dg_p(H\wr G)=\dg_p(H)+\dg_p(G)$$ for any prime $p$.
\end{cor}

So, for groups $A,B,C$ as in Lemma \ref{induction
step}, we have:
$$ \dg_p(A\wr (B \wr C))=\dg_p( (A\times B)\wr C) = \dg_p(A\times B\times C)= \dg_p(A)+\dg_p(B)+\dg_p(C) $$ for every $p$. 
In particular, the epimorphisms in Lemma \ref{induction step}  are $\dg$-preserving.

The same observation holds for Corollary \ref{decending to ascending}, so one has:
\begin{lem} Let $A_1,...,A_r$ be finite abelian groups. Then
$$ \dg_p(A_1\wr (A_2\wr ...\wr A_r))=\dg_p((A_1\wr ... \wr A_{r-1})\wr A_r)=\dg_p(A_1\times... \times A_r)$$ are all $\sum_{i=1}^r\dg_p(A_i)$ for any prime $p$.
\end{lem}
For cyclic groups $A_1,...,A_r$, $\dg_p(A_1 \wr (A_2 \wr ... \wr A_r))$ is simply the number of cyclic groups among $A_1,..,A_r$ whose $p$-part is non-trivial. Thus:
\begin{cor}\label{d of iterated cyclic is max appearence of p} Let $C_1,...,C_r$ be finite cyclic groups and $G=C_1\wr (C_2 \wr ...\wr C_r)$. Then
$\dg(G)=max_{p||{G}|}\dg(C_1(p)\wr (C_2(p) \wr ... \wr C_r(p)))$.
\end{cor}

Let us apply Lemma \ref{induction step} in order to connect between descending iterated wreath products of abelian and cyclic groups:
\begin{prop}\label{connection between abelian and cyclic} Let $A_1,...,A_r$ be finite abelian groups and let $A_i$ have invariant factors $C_{i,j}$ for $j=1,..,l_i$, i.e. $A_i=\prod_{j=1}^{l_i}C_{i,j}$ and $|C_{i,j}|||C_{i,j+1}|$ for any $i=1,...,r$ and $j=1,...,l_i-1$.
Then there is an epimorphism from the descending iterated wreath product $\widetilde{G}:=\wr_{i=1}^r \wr_{j=1}^{l_i} C_{i,j}$ {\rm(here the groups $C_{i,j}$ are ordered lexicographically: $C_{1,1}, C_{1,2},...,C_{1,l_1},C_{2,1},...,C_{r,l_r}$)} to $G:= A_1\wr (A_2 \wr ...\wr A_r)$.
\end{prop}
\begin{proof}

Let us assume $A_1\not=\{0\}$ (otherwise $A_1$ can be simply omitted). Let us prove the assertion by induction on  $\sum_{i=1}^rl_i$. Let $G_2=A_2\wr (A_3 \wr ... \wr A_k)$. Write $A_1=C_{1,1}\times A_1'$. By Lemma \ref{induction step}, there is an  epimorphism $\pi_1:C_{1,1} \wr (A_1'\wr G_2)\ra (C_{1,1}\times A_1')\wr G_2=A_1\wr G_2=G$. 
By applying the induction hypothesis to $A_1',A_2,...,A_r$, there is an epimorphism $\pi_2'$ from the  descending iterated wreath product $\widetilde{G}_2=\wr_{j=2}^{l_1} C_{1,j} \wr (\wr_{i=2}^r \wr_{j=1}^{l_i} C_{i,j})$ to $A_1'\wr G_2$.
By Lemma \ref{functoriality with second argument}, $\pi_2'$ induces an epimorphism $\pi_2:C_{1,1}\wr \widetilde{G}_2\ra C_{1,1} \wr (A_1'\wr G_2)$. Taking the composition $\pi=\pi_2\pi_1$, we obtain the required epimorphism:
$ \pi: \widetilde{G}=C_{1,1}\wr \widetilde{G}_2\ra G. $
\end{proof}
\begin{rem} Note that:
$$ \dg_p(\widetilde{G})=\sum_{i=1}^r\sum_{j=1}^{l_i} \dg_p(C_{i,j})=\sum_{i=1}^r \dg_p(A_i)=\dg_p(G)$$
for every $p$ and hence $\pi$  is $\dg$-preserving.
\end{rem}
Therefore, showing $G$ is a $\dg$-preserving epimorphic image of an iterated wreath product of abelian groups is equivalent to showing $G$ is a $\dg$-preserving epimorphic image of an iterated wreath product of finite cyclic groups.

\section{Wreath length}

%
%

The following lemma is essential for the definition of wreath length:
\begin{lem} Let $G$ be a finite semiabelian group. Then $G$ is a homomorphic
image of a descending iterated wreath product  of finite cyclic
groups, i.e. there are finite cyclic groups $C_1,...,C_r$ and an
epimorphism $C_1\wr (C_2\wr ... \wr C_r)\ra G.$
\end{lem}
\begin{proof} By Proposition \ref{connection between abelian and cyclic} it suffices to show $G$ is an epimorphic image of a descending iterated wreath product of finite abelian groups. We shall prove this claim by induction on $|G|$. The case $G=\{1\}$ is trivial.  By \cite{Den}, $G=A_1H$ with $A_1$ an abelian normal subgroup and $H$ a proper semiabelian subgroup of $G$. First, there is an epimorphism $\pi_1:A_1\wr H \ra A_1H=G$.  By induction there are abelian groups $A_2,..,A_r$ and an epimorphism $\pi_2':A_2 \wr (A_3 \wr ... \wr A_r)\ra H$. By Lemma \ref{functoriality of first argument}, $\pi_2'$ can be extended to an epimorphism $\pi_2:A_1 \wr (A_2\wr ... \wr A_r) \ra A_1 \wr H$. So, by taking the composition $\pi=\pi_1\pi_2$ one obtains the required epimorphism $\pi:A_1 \wr (A_2 \wr ... \wr A_r) \ra G$. 
 \end{proof}
We can now define:
\begin{defn} Let $G$ be a finite semiabelian group.  Define the {\it wreath length} $\wl(G)$ of $G$ to be the smallest positive integer $r$ such that there are finite cyclic groups $C_1,...,C_r$ and an
epimorphism $C_1\wr (C_2 \wr ... \wr C_r) \ra G$.
\end{defn}
Let $\widetilde{G}=C_1\wr (C_2 \wr ... \wr C_r)$ and $\pi:\widetilde{G}\ra G$ an epimorphism. Then by Corollary \ref{d of iterated cyclic is max appearence of p}: $$ \dg(G)\leq \dg(\widetilde{G})\leq r.$$ In particular $\dg(G)\leq \wl(G)$.
\begin{prop}\label{wreath len of wreath prd} Let $C_1,...,C_r$ be nontrivial finite cyclic groups.  Then $\wl(C_1\wr (C_2 \wr ...\wr C_r))=r$.
\end{prop}
Let $\dl(G)$ denote the derived length of a (finite) solvable group $G$, i.e. the smallest positive integer $n$ such that the $n$th higher commutator subgroup of $G$ ($n$th element in the derived series $G=G^{(0)}\geq G^{(1)}=[G,G]\geq \cdots \geq G^{(i)}=[G^{(i-1)},G^{(i-1)}]\geq \cdots$) is trivial.  In order to prove this proposition we will use the following lemma:
\begin{lem}\label{derived length lemma}   Let $C_1,...,C_r$ be nontrivial finite cyclic groups.  Then $\dl(C_1\wr (C_2 \wr ...\wr C_r))=r$.
\end{lem}
\begin{proof}
It is easy (by induction) to see that $\dl(C_1\wr (C_2 \wr ...\wr C_r))\leq r$.  We turn to the reverse inequality.  By Corollary 2.11, it suffices to prove it for the ascending iterated wreath product
$G=(C_1 \wr ... \wr C_{r-1})\wr C_r$. We prove this by induction on $r$.  The case $r=1$ is trivial.  Assume $r\geq 1$. Write $G_1:=(C_1 \wr ... \wr C_{r-2}) \wr C_{r-1}$ so that $G=G_1\wr C_r$. By induction hypothesis, $\dl(G_1)=r-1.$   View $G$ as the semidirect product $G_1^r\rtimes C_r$.  For any $g\in G_1$, the element $t_g:=(g,g^{-1},1,1,...,1)\in G_1^r$ lies in $[G_1^r,C_r]$ and hence in $[G_1^r,C_r]\leq G' \leq G_1^r$.  Let $H=\{t_g|g\in G_1 \}$.  The projection map $G_1^r\rightarrow G_1$ onto the first copy of $G_1$ in $G_1^r$ maps $H$ onto $G_1$.  Since $H\leq G'$, the projection map also maps $G'$ onto $G_1$.  Now $\dl(G_1)=r-1$ by the induction hypothesis.  It follows that $\dl(G')\geq r-1$, whence $\dl(G)\geq r$.
\end{proof}
To prove the proposition, we first observe that  $\wl(C_1\wr (C_2 \wr ...\wr C_r))\leq r$ by definition.  If
$C_1\wr (C_2 \wr ...\wr C_r)$ were a homomorphic image of a shorter descending iterated wreath product
$C_1'\wr (C_2' \wr ...\wr C_s')$, then by Lemma \ref{derived length lemma}, $s=\dl(C_1'\wr (C_2' \wr ...\wr C_s'))\geq \dl(C_1\wr (C_2 \wr ...\wr C_r))=r>s$, contradiction. \qed

Combining Proposition \ref{wreath len of wreath prd} with Corollary \ref{d of iterated cyclic is max appearence of p} we have:
\begin{cor}\label{descriptio of wl=d for iterated} Let $C_1,...,C_r$ be finite cyclic groups and $G=C_1\wr (C_2 \wr ...\wr C_r)$. Then $\wl(G)=\dg(G)$ if and only if there is a prime $p$ for which $p||C_1|,...,|C_r|$.
\end{cor}

We shall now see that all examples of groups $G$ with $\wl(G)=\dg(G)$ arise from Corollary \ref{descriptio of wl=d for iterated}:
\begin{prop}\label{out char of wl} Let $G$ be a finite semiabelian group. Then $\wl(G)=\dg(G)$ if and only if there is a prime $p$, finite cyclic groups $C_1,...,C_r$ for which $p||C_i|$, $i=1,...,r$,  and a $\dg$-preserving epimorphism $\pi:C_1 \wr (C_2 \wr ... \wr C_r)\ra G$.
\end{prop}
\begin{proof} Let $d=\dg(G)$. The equality $d=\wl(G)$ holds if and only if there are finite cyclic groups $C_1, C_2,...,C_d$ and an epimorphism $\pi:\widetilde{G}=C_1 \wr (C_2 \wr ... \wr C_d)\ra G$.  
Assume the latter holds. Clearly  $d\leq \dg(\widetilde{G})$ but by Corollary \ref{d of iterated cyclic is max appearence of p} applied to $\widetilde{G}$ we also have $\dg(\widetilde{G})\leq d$. It follows that $\pi$ is $\dg$-preserving. Since $\dg(G)=max_p(\dg_p(G))$, there is a prime $p$ for which $d=\dg_p(G)$ and hence $\dg_p(\widetilde{G})=d$. Thus, $p||C_i|$ for all $i=1,...,r$.

Let us prove the converse. Assume there is a prime $p$, finite cyclic groups $C_1,...,C_r$ for which $p||C_i|$, $i=1,...,r$,  and a $\dg$-preserving epimorphism $\pi:\widetilde{G}:=C_1 \wr (C_2 \wr ... \wr C_r)\ra G$. Since $p||C_i|$, it follows that $\dg_p(\widetilde{G})=r$. As $\dg_p(\widetilde{G})\leq \dg(\widetilde{G})\leq r$, it follows that $\dg(G)=\dg(\widetilde{G})=r$. In particular $\wl(G)\leq r=\dg(G)$ and hence $\wl(G)=\dg(G)$.

\end{proof}

\begin{rem}\label{cyclic dec of p-groups}
Let $G$ be a semiabelian $p$-group.
By \cite[Corollary 2.15]{N}, $G$ is a $\dg$-preserving image of an iterated wreath product of abelian
subgroups of $G$ (following the proof one can observe that the abelian groups were
actually subgroups of $G$). So, by Proposition \ref{connection between abelian and cyclic},
$G$ is a $\dg$-preserving epimorphic image of $\widetilde{G}:=C_1 \wr (C_2 \wr ... \wr C_k)$
for cyclic subgroups $C_1,...,C_k$ of $G$. By applying Proposition \ref{out char of wl} one obtains $\wl(G)=\dg(G)$.
\end{rem}
\begin{rem} Throughout the proof of \cite[Corollary 2.15]{N} one can use the minimality assumption posed on the decompositions to show directly that the abelian groups $A_1,...,A_r$, for which there is a $\dg$-preserving epimorphism $A_1 \wr (A_2 \wr ... \wr A_r)\ra G$, can be actually chosen to be cyclic.
\end{rem}

We shall generalize Remark \ref{cyclic dec of p-groups} to nilpotent groups:

\begin{prop}\label{wl=d for nilp} Let $G$ be a finite nilpotent semiabelian group. Then $\wl(G)=\dg(G)$.
\end{prop}
\begin{proof} Let $d=\dg(G)$. Let $p_1,...,p_k$ be the primes dividing $|G|$ and let $P_i$ be the $p_i$-Sylow subgroup of $G$ for every $i=1,...,k$. So, $G\cong \prod_{i=1}^k P_i$. By Remark \ref{cyclic dec of p-groups}, there are cyclic $p_i$-groups $C_{i,1},...,C_{i,r_i}$ and a
$\dg$-preserving epimorphism $\pi_i:C_{i,1} \wr (C_{i,2} \wr ... \wr C_{i,r_i})\ra P_i$ for every $i=1,..,k$.
In particular for any $i=1,...,k$, $r_i= \dg(P_i)=\dg_p(G)\leq d$.  For any $i=1,..,k$ and any $d\geq j>r_i$, set $C_{i,j}=\{1\}$. For any $j=1,...,d$
define $C_j=\prod_{i=1}^kC_{i,j}$.

We claim $G$ is an epimorphic image of $\widetilde{G}=C_1 \wr (C_2 \wr ... \wr C_{d})$. To prove this claim it suffices
to show every $P_i$ is an epimorphic image of $\widetilde{G}$ for every $i=1,..,k$. As $C_{i,j}$ is an epimorphic image of $C_j$ for every
$j=1,...,d$ and every $i=1,..,k$, one can apply Lemmas \ref{functoriality of first argument}  and \ref{functoriality with second argument} iteratively to obtain an epimorphism $\pi_i':\widetilde{G}\ra C_{i,1} \wr (C_{i,2} \wr ... \wr C_{i,r})$ for every $i=1,...,k$. Taking the
composition $\pi_i'\pi_i$ gives the required epimorphism and proves the claim.  As $G$ is an epimorphic image of an iterated wreath product of $\dg(G)$ cyclic groups
one has $\wl(G)\leq \dg(G)$ and hence $\wl(G)=\dg(G)$.

\end{proof}

\begin{exam} Let $G=D_n=\langle \sigma,\tau | \sigma^2=1,\tau^n=1, \sigma\tau\sigma=\tau^{-1} \rangle$ for $n\geq 3$. Since $G$ is an epimorphic image of $\langle\tau\rangle\wr \langle\sigma\rangle$ and $G$ is not abelian we have $\wl(G)=2$. 
On the other hand $\dg(G)=\dg(G/[G,G])$ is $1$ if $n$ is odd and $2$ if $n$ is even. So, $G=D_3=S_3$ is the minimal example for which $\wl(G)\not=\dg(G)$.
\end{exam}


\section{a ramification bound for semiabelian groups} In
this section we prove:

\begin{thm}\label{main Theorem} Let $G$ be a
finite semiabelian group.  Then there exists a tamely ramified extension $K/\mQ$ with $G(K/\mQ)\cong G$
in which at most $\wl(G)$ primes ramify.  \end{thm}

The
proof relies on the splitting Lemma from
\cite{KS}: Let $\ell$ be a rational prime, $K$ a number field and
$\fp$ 
a prime of $K$ that is prime to $\ell$.
Let $I_{K,\fp}$ denote the group of fractional ideals prime to
$\fp$, $P_{K,\fp}$ the subgroup of principal ideals that are prime
to $\fp$ and let  $P_{K,\fp,1}$ be the subgroup of principal ideals
$(\alpha)$ with $\alpha\equiv 1$ (mod $\fp$). Let $\overline
P_{\fp}$ denote $P_{K,\fp}/P_{K,\fp,1}$. The ray class group
$Cl_{K,\fp}$ is defined to be $I_{K,\fp}/P_{K,\fp,1}$. Now, as
$I_{K,\fp}/P_{K,\fp}\cong Cl_{K}$, one has the following short exact
sequence: \begin{equation} \label{ray class field sequence}
1 \longrightarrow \overline
P_{\fp}^{(\ell)} \longrightarrow Cl_{K,\fp}^{(\ell)} \longrightarrow
Cl_{K}^{(\ell)} \longrightarrow 1, \end{equation} where $A^{(\ell)}$
denotes the $\ell$-primary component of an abelian group $A$. Let us
describe a sufficient condition for the splitting of (\ref{ray class
field sequence}). Let $\fa_1,...,\fa_r \in I_{K,\fp}$,
$\tilde\fa_1,...,\tilde\fa_r$ their classes in $Cl_{K,\fp}^{(\ell)}$
with images $\overline{\fa}_1,...,\overline{\fa}_r$ in
$Cl_K^{(\ell)}$, so that
$Cl_{K}^{(\ell)}=\langle\overline{\fa}_1\rangle\times
\langle\overline{\fa}_2\rangle\times...\times \langle
\overline{\fa}_r\rangle$.  Let $\ell^{m_i}:=|\langle \overline
\fa_i\rangle |$ and let $a_i\in K$ satisfy
$\fa_i^{\ell^{m_i}}=(a_i)$, for $i=1,...,r$.
\begin{lem}\label{useful part of the splitting
lemma}\emph{(Kisilevsky-Sonn \cite{KS2})} Let $\fp$ be a prime of
$K$ and let $K'=K(\sqrt[\ell^{m_i}]{a_i}|i=1,...,r)$. If $\fp$
splits completely in $K'$ then the sequence (\ref{ray class field
sequence}) splits. \end{lem} The splitting  of (\ref{ray class field
sequence}) was used in \cite{KS} to construct cyclic ramified
extensions at one prime only. Let $m=\max\{1,m_1,...m_r\}$.  Let
$U_K$ denote the units in $\cO_K$. \begin{lem}\label{cor-existence
of totally ramified extension}\emph{(Kisilevsky-Sonn \cite{KS})} Let
$K''=K(\mu_{\ell^m},\sqrt[\ell^m]{\xi},\sqrt[\ell^{m_i}]{a_i}|\
\xi\in U_K, i=1,...,r)$ and $\fp$ a prime of $K$ which splits
completely in $K''$. Then there is a cyclic $\ell^m$-extension of
$K$ that is totally ramified at $\fp$ and is not ramified at any
other prime of $K$. \end{lem}

\begin{cor}\label{K'''} Let $K$ be a number field, $n$ a positive integer.  Then there exists a finite extension $K'''$ of $K$ such that if $\fp$ is any prime of $K$ that splits completely in $K'''$, then there exists a cyclic extension $L/K$ of degree $n$ in which $\fp$ is totally ramified and $\fp$ is the only prime of $K$ that ramifies in $L$.
\end{cor}
\begin{proof} Let $n=\prod_{\ell}\ell^{m(\ell)}$ be the decomposition of $n$ into primes.  Let $K'''$ be the composite of the fields $K''=K''(\ell)$ in Lemma \ref{cor-existence of totally ramified extension} ($m=m(\ell))$.  Let $L(\ell)$ be the cyclic extension of degree $\ell^{m(\ell)}$ yielded by Lemma \ref{cor-existence of totally ramified extension}.  The composite $L=\prod L(\ell)$ has the desired property.
\end{proof}


\begin{proof}(Theorem \ref{main Theorem})  \rm By definition, $G$ is a homomorphic image of a descending iterated wreath product of cyclic groups $C_1 \wr (C_2 \wr \cdots \wr C_r)$, $r=\wl(G)$.  Without loss of generality $G\cong C_1 \wr (C_2 \wr \cdots \wr C_r)$ is itself a  descending iterated wreath product of cyclic groups.  Proceed by induction on $r$.  For $r=1$, $G$ is cyclic of order say $N$. If $p$ is a rational prime $\equiv 1$
(mod $N$), then the field of $p$th roots of unity $\mQ(\mu_p)$
contains a subfield $L$ cyclic over $\mQ$ with Galois group $G$ and
exactly one ramified prime, namely $p$.  Thus the theorem holds for
$r=1$.

Assume $r>1$ and the theorem holds for $r-1$.  Let $K_1/\mQ$ be a tamely ramified Galois extension with $G(K_1/\mQ)\cong G_1$, where $G_1$ is the descending iterated wreath product $C_2 \wr (C_3 \wr \cdots \wr C_r)$, such that the ramified primes in $K_1$ are a subset of $\{p_2,...,p_r\}$.  By Corollary \ref{K'''}, there exists a prime $p=p_1$ not dividing the order of $G$ which splits completely in $K_1'''$, the field supplied for $K_1$ by Corollary \ref{K'''}, and let $\fp=\fp_1$ be a prime of $K_1$ dividing $p$.  By Corollary \ref{K'''}, there exists a cyclic extension $L/K_1$ with $G(L/K_1)\cong C_1$ in which $\fp$ is totally ramified and in which $\fp$ is the only prime of $K_1$ which ramifies in $L$.

Now $\fp$ has $|G_1|$ distinct conjugates $\{\sigma(\fp)|\sigma \in
G(K_1/\mQ)\}$ over $K_1$.  For each $\sigma \in G(K_1/\mQ)$, the conjugate extension
$\sigma(L)/K_1$ is well-defined, since $K_1/\mQ$ is Galois.  Let $M$ be the
composite of the   $\sigma(L)$, $\sigma \in G(K_1/\mQ)$.  For each $\sigma$,
$\sigma(L)/K_1$ is cyclic of degree $|C_1|$, ramified only at $\sigma(\fp)$, and
$\sigma(\fp)$ is totally ramified in $\sigma(L)/K_1$.  It now follows (see e.g. \cite[Lemma 1]{KS}) that the fields $\{\sigma(L)|\sigma \in G(K_1/\mQ)\}$ are linearly disjoint over
$K_1$, hence $G(M/\mQ)\cong C_1 \wr  G_1 \cong G$. Since the only primes of
$K_1$ ramified in $M$ are $\{\sigma(\fp)|\sigma \in G(K_1/\mQ)\}$, the only rational
primes ramified in $M$ are $p_1,p_2,...,p_n$. \end{proof}

\begin{cor} The minimal ramification problem has a positive solution for all finite semiabelian groups $G$ for which $\wl(G)=\dg(G)$.  Precisely, any finite semiabelian group $G$ for which $\wl(G)=\dg(G)$ can be realized tamely as a Galois group over the rational numbers with exactly $\dg(G)$ ramified primes. \end{cor}

By Proposition \ref{wl=d for nilp}, we have

\begin{cor} The minimal ramification problem has a positive solution for all finite nilpotent semiabelian groups.  \end{cor}

\section{Arithmetic consequences}

\bigskip
In this section we examine some arithmetic consequences of a
positive solution to the minimal ramification problem.
Specifically, given a group $G,$ the existence of infinitely many
minimally tamely ramified $G$-extensions $K/\Q$ is re-interpreted
in some cases in terms of some open problems in algebraic
number theory. We will be most interested in the case $\dg(G)=1.$

\begin{prop}{\label{prop:classnumber divisibility}
Let $q$ and $\ell$ be distinct primes. Let $K/\Q$ be a cyclic extension of degree  $n:=[K:\Q]\geq 2$ with $(n,q\ell)=1$. Suppose that $K/\Q$ is totally and tamely ramified  at a unique prime $\mathfrak l$ dividing $\ell$.  Then $q$ divides the class number $h_K$ of $K$ if and only if
there exists an extension $L/K$ satisfying the following:
\item{i).} $L/\Q$ is a Galois extension with {\bf non-abelian} Galois group $G=G(L/\Q)$.
\item{ii).} The degree $[L:K]=q^s$ is a power of $q.$
\item{iii).} $L/\Q$ is (tamely) ramified only at primes over $\ell.$}
\end{prop}

\begin{proof}
 First suppose that $q$ divides $h_K.$  Let $K_0$ be the $q$-Hilbert class field
of $K,$ {\it i.e.\/} $K_0/K$ is the maximal unramified abelian $q$-extension of $K.$
Then $K_0/\Q$ is a Galois extension with Galois group $G:=G(K_0/\Q),$
and $H:=G(K_0/K)\simeq (C_K)_q \neq 0$, the $q$-part
of the ideal class group of $K.$ Then $[G,G]$ is contained in $H.$
If $[G,G] \subsetneq H,$ then the fixed field of $[G,G]$ would be an abelian extension of $\Q$ which
contains an unramified $q$-extension of $\Q$ which is impossible. Hence $[G,G]=H\neq 0$ and so $G$ is
a non-abelian group, and $L=K_0$ satisfies {\it i),ii)}, and {\it iii)} of the statement.

Conversely suppose that there is an extension $L/K$ satisfying {\it i),ii)}, and {\it iii)} of the statement.
Since $H=G(L/K)$ is a $q$-group, there is a sequence of normal subgroups
$H=H_0\supset H_1\supset H_2 \cdots \supset H_s=0$ with $H_i/H_{i+1}$ a cyclic
group of order $q.$ Let $L_i$ denote the fixed field of $H_i$ so that
$K=L_0\subset \cdots\subset L_s=L.$ Let $m$ be the largest index such that $L_m/\mQ$ is totally ramified
(necessarily at $\ell$). If $m=s,$ then $L/\Q$ is totally and tamely ramified at $\ell$ and so the inertia
group $T(\mathfrak L/(\ell))=G,$ where in this case $\mathfrak L$ is the unique prime of $L$ dividing $\ell.$
Since $L/\Q$ is tamely ramified it follows that  $T(\mathfrak L/(\ell))$ is cyclic, but this contradicts
the hypothesis that $G$ is non-abelian. Therefore it follows that $m<s,$ and so $L_{m+1}/L_m$ is unramified
and therefore $q$ must divide the class number $h_{L_m}.$ Then a result of Iwasawa \cite{I} implies that $q$
divides all of the class numbers $h_{L_{m-1}}, \cdots,h_{L_0}=h_K.$
\end{proof}

We now apply this to the case that $G\not=\{1\}$ is a quotient of the regular wreath product $C_q \wr C_p$
where $p$ and $q$ are distinct primes. Then $\dg(G)=1.$

The existence of infinitely many minimally tamely ramified $G$-extensions $L/\Q$ would by
Proposition \ref{prop:classnumber divisibility} imply the existence of infinitely many
cyclic extensions $K/\Q$ of degree $[K:\Q]=p$ ramified at a unique prime $\ell \neq p,q$
for which $q$ divides the class number $h_K.$  
(If there were only finitely many distinct such cyclic extensions $K/\Q,$
then the number of ramified primes $\ell$ would be bounded, and there would be an absolute upper bound
on the possible discriminants of the distinct fields $L/\Q.$ By Hermite's theorem, this would mean that
the number of such $G$-extensions $L/\Q$ would be bounded).

The question of whether there is an infinite number of cyclic degree $p$ extensions (or even one) of $\mQ$ whose
class number is divisible by $q$ is in general open at this time.

For $p=2,$ it is known that
there are infinitely many quadratic fields (see Ankeny, Chowla \cite{AC}), with class numbers divisible
by $q,$ but it is not known that this occurs for quadratic fields with prime discriminant.

This latter statement is also a consequence of Schinzel's hypothesis as is shown by Plans in \cite{P}.
There is also some numerical evidence that the heuristic of Cohen-Lenstra should be statistically independent of
the primality of the discriminant (see Jacobson, Lukes, Williams \cite{JLW} or te Riele, Williams \cite{RW}).
If this were true, then one would expect that there is a positive density of primes $\ell$
for which the cyclic extension of degree $p$ and conductor $\ell$ would have class number divisible by $q.$

For $p=3$ it has been proved by Bhargava \cite{B} that there are infinitely many cubic fields $K/\Q$
for which $2$ divides their class numbers. That there are infinitely many cyclic cubics with prime
squared discriminants whose class numbers are even (or more generally divisible by some fixed prime $q$)
seems out of reach at this time.

In our view, there is a significant arithmetic interest in solving the minimal ramification problem for other groups (see also \cite{H}, \cite{JR}, \cite{R}).

\end{document}